\documentclass[12pt]{article}

\usepackage[english]{babel}

\usepackage[letterpaper,top=2cm,bottom=2cm,left=3cm,right=3cm,marginparwidth=1.75cm]{geometry}

\usepackage{amsmath}
\usepackage{amssymb}
\usepackage{rotating}
\usepackage{mathrsfs}
\usepackage{graphicx}
\usepackage[colorlinks=true, allcolors=blue]{hyperref}
\usepackage{natbib}
\usepackage{mdframed}
\usepackage{url}
\usepackage{array}

\def\thickhline{\noalign{\hrule height 1.5pt}}

\title{AdHImEx: Adaptively High-Order Implicit-Explicit Transport for Large Time Steps}
\author{Amber J. te Winkel, Hilary Weller, Christian K\"uhnlein, James Kent}
\date{September 8, 2026}

\begin{document}
\maketitle
\newcolumntype{?}{!{\vrule width 1.5pt}}

\section*{Abstract}

Adaptively Implicit-Explicit (AdImEx) time stepping provides stability for large time steps for mass-conservative transport, making it attractive for the numerical representation of advection in weather and climate prediction. Existing AdImEx transport schemes become first-order accurate in the large-time-step, implicit regime and may require multiple sparse matrix solutions per time step. This work introduces AdHImEx, a new AdImEx scheme that addresses both limitations simultaneously. AdHImEx provides second-order accuracy for large Courant numbers while requiring only a single matrix solution per time step. It is a Runge-Kutta scheme that blends third-order accurate explicit time stepping and Crank-Nicolson implicit time stepping to provide numerically verified stability for Courant numbers up to 100. AdHImEx is designed for flows that are predominantly explicit, activating implicit time stepping only where large Courant numbers occur. It retains the full efficiency and third-order accuracy of the explicit time stepping in regions with small Courant numbers. In this work, it is combined with a fifth-order accurate finite-volume discretisation in space for mass conservation. A second contribution of this work is a novel stage-dependent treatment of the divergence operator that preserves constancy despite different implicit and explicit Runge-Kutta stage time steps and spatially varying implicitness, removing spurious divergence. This innovation is applicable to other AdImEx time stepping schemes. Convergence, transport, and efficiency tests demonstrate substantial improvements over previous first-order AdImEx schemes, little phase and amplitude error, accurate transport on locally refined meshes with large Courant numbers, and promising efficiency characteristics for atmospheric transport.

\section{Introduction}

To solve the dynamical equations governing atmospheric fluid motion, numerical weather prediction and climate models need to represent the linear advection equation,
\begin{equation}
    \frac{\partial \Psi}{\partial t} + \nabla\cdot (\vec{u} \Psi) = 0
    \label{eqadv}
\end{equation}
for a scalar field $\Psi$ in space $\vec{x}$ and time $t$ with velocity $\vec{u}$. $\Psi$ is a density, either dry density $\rho$ or the tracer mass density $\rho q$ (the tracer concentration $\rho$ weighted by the dry density). Important properties for numerical advection schemes include stability, mass conservation, accuracy, and efficiency \citep{mengaldoCurrentEmergingTimeIntegration2019}.
It is not straightforward to incorporate all of these at the same time, in particular on unstructured grids.

Local mass conservation can be guaranteed by using a finite-volume discretisation, which calculates the change in $\Psi$ over time by considering the fluxes going in and out of grid cells. Semi-Lagrangian time stepping can be combined with finite-volume in space for a flux-form semi-Lagrangian scheme (FFSL). FFSL offers favourable stability, restricted by the Lipschitz rather than Courant number, and can be accurate and efficient. Multi-dimensional FFSL schemes were presented by \citet{linMultidimensionalFluxFormSemiLagrangian1996} and \citet{leonardConservativeExplicitUnrestrictedTimeStep1996}. \citet{bendallSWIFTMonotonicFluxForm2025} built on this, introducing multi-dimensional monotonicity, nonnegativity, and constancy preservation for density-weighted fields. However, even though FFSL schemes can be accurate and efficient on structured grids, they can become complicated and computationally expensive on arbitrary, unstructured grids. This was addressed by \citet{herzfeldFluxFormSemiLagrangianAdvection2023}, who adapted FFSL to work on arbitrary grids but note that optimised parallel integration remains an issue.

To cater to arbitrary grids, one could use explicit Eulerian time stepping with finite-volume in space. However, the favourable stability of semi-Lagrangian is lost: Where semi-Lagrangian time stepping remains stable for large time steps, (explicit) Eulerian may go unstable \citep{durranNumericalMethodsFluid2010}. This is problematic as efficiency asks for as long a time step as possible. For Eulerian time stepping, stability is indicated by the Courant number $C$, describing the number of cells travelled in one time step and depending on the velocity $\vec{u}$, time step $\Delta t$, and cell volume $V$. Increasing $\Delta t$ corresponds to increasing the Courant number and may lead to instabilities.

To calculate the tracer field $\Psi$ at the time level $(n+1)$, explicit Eulerian schemes use $\Psi$ at time level $(n)$ (and earlier for multistep schemes). They require the Courant number to be restricted to, e.g., one for stability. Implicit schemes also use the field at $(n+1)$ and may be unconditionally stable, i.e., they do not have a maximum Courant number for stability that restricts the time step size. They do require the solution of a matrix, reducing efficiency per time step compared to explicit time stepping \citep{mengaldoCurrentEmergingTimeIntegration2019}.

One way around the instability problem of explicit Eulerian time stepping is to use adaptively implicit-explicit (AdImEx) time stepping \citep[e.g.,][]{wellerAdaptivelyImplicitMPDATA2023}. AdImEx time stepping combines explicit and implicit Eulerian time stepping such that the result is stable by using implicit time stepping where needed. Whenever the explicit time stepping would go unstable, it locally blends with an unconditionally stable implicit time stepping to ensure stability. Setting the time step such that the explicit stability limit is rarely exceeded, AdImEx schemes could be able to largely retain the higher order of accuracy and efficiency of explicit time stepping. When AdImEx time stepping is combined with a finite-volume discretisation in space, the result can be conservative, stable, accurate, and work on arbitrary grids. By removing the time step restriction, AdImEx time stepping may give overall efficiency gains compared to fully explicit time stepping when paired with an appropriate sparse matrix solver. This is particularly attractive for variable-resolution meshes, where a small number of highly refined cells would otherwise (with an explicit scheme) impose a restrictive global time step.

Recent AdImEx schemes were first-order accurate for large time steps in time, space, or both. \citet[WS20]{wickerImplicitExplicitVertical2020}, \citet{liEnhancingStabilityGlobal2022}, and \citet{liLongTimeStepPermittingTracerTransport2024} employed AdImEx time stepping in the vertical to create stable and robust schemes for large time steps in the presence of strong vertical velocities. WS20 showed improved efficiency for convection-permitting resolutions. \citet{wellerAdaptivelyImplicitMPDATA2023,wellerAdaptivelyImplicitAdvection2024} looked at two-dimensional AdImEx tracer transport on arbitrary grids for locally very large Courant numbers. However, all these only used first-order accurate implicit time stepping. Furthermore, \citet{wellerAdaptivelyImplicitAdvection2024} required multiple matrix solutions each time step, inhibiting efficiency.

Most recently, the AdLImHEx scheme was presented \citep{wellerAdaptivelyImplicitExplicitDensityWeightedsubmitted}. This scheme blends fourth-order in time and fifth-order in space for small Courant numbers with first-order in time as Courant numbers grow large. Importantly, it requires one solve per Runge-Kutta stage (of which there are four), and was demonstrated on the sphere with multiple grid types and density-weighted tracers. The result was remarkably efficient for large Courant numbers. However, \citet{wellerAdaptivelyImplicitExplicitDensityWeightedsubmitted} still has the above two limitations: first-order accuracy for the implicit part and the use of multiple implicit solves per time step.

This paper makes two methodological contributions. First, it introduces AdHImEx, a new AdImEx transport scheme that combines second-order temporal accuracy for large Courant numbers with only a single implicit solve per time step. To our knowledge, this has not been done before. The time stepping is combined with a finite-volume discretisation, giving a conservative scheme with stability for large Courant numbers, third- or second-order accuracy in time, and fifth-order accuracy in space.

Second, it derives a stage-dependent treatment of the divergence operator that preserves constancy. AdHImEx is required to take Runge-Kutta stage time steps with different sizes for its implicit and explicit parts. This introduces a challenge because spatially varying implicitness introduces spurious divergence, destroying constancy preservation. This work derives a stage-dependent treatment of the divergence operator that preserves both conservation and constancy. Importantly, this innovation is not specific to AdHImEx. It provides a framework for preserving constant fields in a broader class of AdImEx schemes, increasing the flexibility in which implicit and explicit schemes can be combined.

One-dimensional results first verify the expected accuracy properties of AdHImEx and quantify improvements over previous AdImEx time stepping. Two-dimensional tests then demonstrate the practical impact of the method, including transport on locally refined meshes with large Courant numbers, constancy preservation, and efficiency. Finally, conclusions and future work are discussed.

\section{Numerical Methods}\label{secmethods}

\subsection{AdHImEx Time Stepping}\label{timestepping}

AdHImEx time stepping is Adaptively Implicit-Explicit time stepping with High order accuracy for both the explicit and implicit parts. It combines two Runge-Kutta schemes, one explicit and one diagonally implicit. Runge-Kutta schemes are multi-stage schemes, which, at each time step, calculate multiple intermediate fields $\Psi^s$, called `stages', that are then combined to find $\Psi^{(n+1)}$. This can be written as
\begin{equation}
    \begin{aligned}
        \Psi^s       & = \Psi^{(n)} - \Delta t \sum^{s}_{r=1} \mathscr{F}_s \Big( a_{sr}, \vec{u}^r, \Psi^r \Big), \qquad s=1,...,\mathfrak{S} \\
        \Psi^{(n+1)} & = \Psi^{(n)} - \Delta t \sum^{\mathfrak{S}}_{s=1} \mathscr{F} \Big( b_{s}, \vec{u}^s, \Psi^s \Big).
    \end{aligned}
    \label{eqrungekutta}
\end{equation}
$\mathscr{F}_s$ and $\mathscr{F}$ are functions, for linear advection specifying the spatial discretisation of the divergence. $\mathfrak{S}$ is the total number of stages, and $(n)$ and $\Delta t$ identify the time level and time step size, respectively. Specific Runge-Kutta time stepping methods are defined by the $a_{sr}$ and $b_s$ coefficients in (\ref{eqrungekutta}). Their values are given in a Butcher tableau,
\begin{equation}
    \begin{array}
        {c|ccccc}
        c_1               & a_{11}               & 0                    & \cdots & 0                               \\
        c_2               & a_{21}               & a_{22}               & \cdots & 0                               \\
        \vdots            & \vdots               & \vdots               & \ddots & \vdots                          \\
        c_{_\mathfrak{S}} & a_{_{\mathfrak{S}1}} & a_{_{\mathfrak{S}2}} & \cdots & a_{_{\mathfrak{S}\mathfrak{S}}} \\
        \hline
                          & b_1                  &
        b_2               & \cdots               & b_{_\mathfrak{S}}
    \end{array}
    \label{butchertableau}
\end{equation}
The number of $a$ rows is the number of stages, and the $b$ row provides the weights for each stage field at the end of the time step, essentially denoting an inherently explicit stage itself. The diagonal $a_{ss}$ coefficients are connected to the implicit terms as can be seen from (\ref{eqrungekutta}). For explicit time stepping, all of these are zero. For diagonally implicit time stepping, at least one will be nonzero. $c_s = \sum^{\mathfrak{S}}_{r=1} a_{sr}$ denotes the time step of the stage: $\Psi^{s}$ is at time $t^{(n)} + c_s \Delta t$. (\ref{eqrungekutta}) does not expressly use this, as the winds $\vec{u}^r=\vec{u}$ are assumed constant in a single time step. In particular, we use the wind velocity at the mid-point $\vec{u}=\vec{u}^{(n+1/2)}$ to calculate $\Psi^{(n+1)}$. This is a centred-in-time, second-order accurate approximation, and is required to preserve constant fields (Section \ref{secconstancy}).

AdHImEx combines the two-stage, second-order, unconditionally stable, implicit Crank-Nicolson time stepping \citep{crankPracticalMethodNumerical1947}
\begin{equation}
    \begin{array}
        {c|ccccc}
        0 & 0           & 0           \\
        1 & \frac{1}{2} & \frac{1}{2} \\
        \hline
          & \frac{1}{2} & \frac{1}{2} \\
    \end{array}
    \label{eqbutcherCN}
\end{equation}
with high-order, explicit Runge-Kutta time stepping. In particular, we use the three-stage explicit scheme from \citet{gottliebStrongStabilityPreservingHighOrder2001} with Butcher tableau,
\begin{equation}
    \begin{array}
        {c|ccccc}
        0           & 0           & 0           & 0           \\
        1           & 1           & 0           & 0           \\
        \frac{1}{2} & \frac{1}{4} & \frac{1}{4} & 0           \\
        \hline
                    & \frac{1}{6} & \frac{1}{6} & \frac{2}{3} \\
    \end{array}
    \label{eqbutcherRK3}
\end{equation}
which is Strong Stability Preserving, or Total Variation Diminishing, for Courant numbers less than or equal to 1. These two tableaus are combined with a blending factor $\theta\in[0,1]$ representing the implicitness:
\begin{equation}
    \begin{array}
        {c|ccccc}
        0                     & 0                 & 0                     & 0                     & 0                     & 0                 \\
        \frac{1}{2}\theta     & \frac{1}{2}\theta & 0                     & 0                     & 0                     & 0                 \\
        1 - \frac{1}{2}\theta & \frac{1}{2}\theta & 1-\theta              & 0                     & 0                     & 0                 \\
        \frac{1}{2}           & \frac{1}{2}\theta & \frac{1}{4}(1-\theta) & \frac{1}{4}(1-\theta) & 0                     & 0                 \\
        1                     & \frac{1}{2}\theta & \frac{1}{6}(1-\theta) & \frac{1}{6}(1-\theta) & \frac{2}{3}(1-\theta) & \frac{1}{2}\theta \\
        \hline
                              & \frac{1}{2}\theta & \frac{1}{6}(1-\theta) & \frac{1}{6}(1-\theta) & \frac{2}{3}(1-\theta) & \frac{1}{2}\theta
    \end{array}
    \label{eqbutcherAdHImEx}
\end{equation}
The specific combination of explicit and implicit terms mimics the approach from \citet{ullrichOperatorSplitRungeKutta2012}, who apply Strang splitting to their explicit and implicit operators. Their scheme was written as a Butcher tableau by \citet{wellerRungeKuttaIMEX2013}.

For $\theta=0$, (\ref{eqbutcherAdHImEx}) reverts to (\ref{eqbutcherRK3}) and the time stepping is fully explicit. For $\theta=1$, it reverts to (\ref{eqbutcherCN}) and the time stepping is fully implicit, with only one implicit solve per time step, in the last stage. A blending factor $\theta$ will be found (Section \ref{stability}) that depends on the Courant number, varies in space and time, and gives AdHImEx stability for large Courant numbers.

\subsection{Constancy Preservation}\label{secconstancy}

For a divergence-free velocity field, a tracer that is initially uniform should remain uniform. Following \citet{leonardConservativeExplicitUnrestrictedTimeStep1996}, we refer to this property as constancy preservation. Constancy preservation serves to prevent the generation of spurious tracer features arising solely from the discretisation of the divergence term. It addresses the same underlying issue as the consistency-with-continuity property discussed by \citet{wellerAdaptivelyImplicitExplicitDensityWeightedsubmitted}, whether this corresponds to a spatially uniform tracer remaining constant or to a conserved variable evolving consistently with the density field. To have either, the velocity used has to remain constant within each time step. This applies to both AdHImEx and AdLImHEx \citep{wellerAdaptivelyImplicitExplicitDensityWeightedsubmitted}, and both use the mid-time velocity for all divergence operators in a time step, giving second-order accuracy.

The AdHImEx scheme, however, has an additional requirement. The Runge-Kutta coefficients for stages $s=2$ and $s=3$ depend on the local value of $\theta$ through (\ref{eqbutcherAdHImEx}): The explicit and implicit parts take different stage time steps. Since $\theta$ may vary spatially, these coefficients are generally not spatially uniform. When such spatially varying coefficients appear inside a flux-form divergence operator, they can induce a spurious divergence even when the physical velocity field itself is divergence-free. An initially uniform tracer field may no longer remain uniform. To ensure both conservation and constancy preservation, a careful choice must therefore be made for the operators $\mathscr{F}_s$ and $\mathscr{F}$.

For conservation, $\mathscr{F}$ needs to use the divergence operator in flux form,
\begin{equation}
    \mathscr{F} \Big( b_s, \vec{u}, \Psi^s \Big)         = \nabla \cdot \Big( b_s \vec{u} \Psi^s \Big).
    \label{eqFfinal}
\end{equation}
However, if $\mathscr{F}_{s}$ is in flux form for all stages, i.e.,
\begin{equation}
    \mathscr{F}_{s} \Big( a_{sr}, \vec{u}, \Psi^r \Big) = \nabla \cdot \Big( a_{sr} \vec{u} \Psi^r \Big) \text{ for } s=1,2,3,4,5,
    \label{eqFs_ncp}
\end{equation}
then AdHImEx does not preserve constancy. Why is this the case? Constancy is preserved when an initially uniform field remains uniform under nondivergent flows \citep{leonardConservativeExplicitUnrestrictedTimeStep1996}. At each grid point, zero rate of change with time requires both
\begin{align}
    \sum^{s}_{r=1} \mathscr{F}_{s} \Big( a_{sr}, \vec{u}, \Psi^r \Big)     & \equiv 0 \quad \forall s, \quad \text{and} \\
    \sum^{\mathfrak{S}}_{s=1} \mathscr{F} \Big( b_s, \vec{u}, \Psi^s \Big) & \equiv 0.
\end{align}
with $\Psi^{(n)}=\Psi^{r=1}$ constant in space, $\nabla \cdot \vec{u} = 0$, and $\vec{u}$ constant within each time step. For general $\gamma$, $\Gamma$, $\mathscr{F}$, $\alpha_\gamma$, and $\Psi$, we need
\begin{equation}
    \sum^{\Gamma}_{\gamma=1} \mathscr{F} \Big( \alpha_\gamma, \vec{u}, \Psi \Big) \equiv 0.
\end{equation}
For the flux-form divergence for $\mathscr{F}$, this reduces to,
\begin{equation}
    \sum^{\Gamma}_{\gamma=1} \nabla \cdot \Big( \alpha_\gamma \vec{u} \Psi \Big) = \sum^{\Gamma}_{\gamma=1} \alpha_\gamma \Psi \nabla \cdot \vec{u} + \sum^{\Gamma}_{\gamma=1} \vec{u} \cdot \nabla \Big( \alpha_\gamma \Psi \Big) = \Psi \vec{u}  \cdot \nabla \Big( \sum^{\Gamma}_{\gamma=1} \alpha_\gamma \Big) \equiv 0.
\end{equation}
As such, we need $\nabla ( \sum^{\Gamma}_{\gamma=1} \alpha_\gamma) \equiv 0$. The coefficient sum ($\widehat{=}\, c_s$ in (\ref{butchertableau})) needs to be spatially uniform. For AdHImEx, this is the case for $b$, and for $a$ when $s=1,4,$ or 5. For $s=2,3$, the dependency on $\theta$ creates spurious divergence. That is, $\nabla \cdot (a_{sr} \vec{u}) \neq 0$ when $\nabla \cdot \vec{u} = 0$. To avoid the resulting spurious divergence, $\alpha_\gamma$ needs to be taken outside of the divergence operator. Then,
\begin{equation}
    \sum^{\Gamma}_{\gamma=1} \alpha_\gamma \nabla \cdot \Big( \vec{u} \Psi \Big) = \sum^{\Gamma}_{\gamma=1} \alpha_\gamma \Psi \nabla \cdot \vec{u} + \sum^{\Gamma}_{\gamma=1} \alpha_\gamma \vec{u} \cdot \nabla \Psi \equiv 0,
\end{equation}
which is the case, preserving constancy. $\mathscr{F}_{s}$ should thus depend on the stage:
\begin{equation}
    \mathscr{F}_{s} \Big( a_{sr}, \vec{u}, \Psi^r \Big) = \begin{cases}
        \nabla \cdot \Big( a_{sr} \vec{u} \Psi^r \Big) & \text{ for } s=1,4,5, \\
        a_{sr} \nabla \cdot \Big( \vec{u} \Psi^r \Big) & \text{ for } s=2,3.
        \label{eqFs}
    \end{cases}
\end{equation}
Crucially, using (\ref{eqFs}) does not break local mass conservation. Even though $\Psi^{s=2}$ and $\Psi^{s=3}$ are no longer conservative, the final result still is. $\mathscr{F}$ has all three variables inside of the divergence, and conservation is maintained where it matters.

We believe this is the first treatment of constancy preservation for spatially varying AdImEx Runge-Kutta coefficients. The impact of the proposed constancy-preserving formulation is demonstrated both through removal of transport artefacts (Section \ref{secartefacts}) and exact preservation of uniform fields (Section \ref{res2dconstancy}).

The same approach is applicable to a broader class of AdImEx time-stepping schemes for preserving constancy. For AdHImEx, we could have chosen to have taken $a_{sr}$ outside of the divergence for $s=1,4,$ and 5 as well. The sole constraint is that $b_s$ (and thus $\theta$) are inside of the divergence when forming $\Psi^{(n+1)}$ to ensure conservation. However, note that this can always be done without problems: $\sum_{s=1}^\mathfrak{S} b_s$ is always independent of $\theta$.

A similar problem is addressed in \citet{melvinMixedFiniteelementFinitevolume2024} for constancy preservation. In this case, constancy preservation is lost through dimensional splitting, introducing spurious divergence for nondivergent flow. In their advective-then-flux approach, this is avoided by writing only the last divergence operator in flux form (for conservation), with all previous operators in advective form. The present work does not consider dimensional splitting and chooses the flux-form operator whenever possible.

\subsection{Spatial Discretisation}\label{spatdisc}

AdHImEx is a method-of-lines scheme, meaning that space and time are discretised independently. We are using a fifth-order finite-volume discretisation in space. To obtain this, we look for an approximation for $\Psi_f$, the tracer field at cell faces. In line with the finite-volume framework, the divergence terms from (\ref{eqFfinal}) and (\ref{eqFs}) become
\begin{equation}
    \nabla \cdot \Big( \vec{\nu} \Psi \Big)_c = \frac{1}{V_c} \sum_{f\in c}\vec{S}_f \cdot \vec{\nu}_f \Psi_f
    \label{eqdiv}
\end{equation}
for a cell $c$. The sum is over the faces $f$ of cell $c$. $\vec{S}_f$ is the outward-pointing face area normal vector, and $\vec{\nu}_f$ is either $a_{sr,f}\vec{u}_f$ or $\vec{u}_f$, assumed to be readily defined at faces. $\Psi_f$, on the other hand, needs to be derived from its cell-centred values. This work uses the same accurate spatial discretisation as for AdLImHEx's high-order explicit part \citep{wellerAdaptivelyImplicitExplicitDensityWeightedsubmitted}. The cell-centred discretisation is interpolated to faces such that when combined with a spatially uniform $\vec{\nu}$, the result is a fifth-order accurate approximation of the divergence in (\ref{eqdiv}) \citep{trembackForwardTimeUpstream1987}. It is given by
\begin{equation}
    \Psi_f = -\tfrac{1}{20}\Psi_{\text{d}^2} + \tfrac{9}{20}\Psi_{\text{d}} + \tfrac{47}{60}\Psi_{\text{u}} - \tfrac{13}{60}\Psi_{\text{u}^2} + \tfrac{1}{30}\Psi_{\text{u}^3},
    \label{eqSD}
\end{equation}
where $\text{d}^*$ and $\text{u}^*$ refer to the $*$th downwind and upwind cells of the face. This is an upwind-biased $[\text{u}^3, \text{d}^2]$ discretisation. We have chosen this particular spatial discretisation because of its good stability properties with the AdHImEx time stepping.

Contrary to AdLImHEx, AdHImEx uses this discretisation not only for the explicit but also for the implicit part. This increases the cost of the matrix solution, but importantly, allows us to assign any errors that arise to the time discretisation - useful for a proof-of-concept work on a new time stepping method. As the spatial discretisation is greater than third-order accurate, errors are dominated by time rather than space.

For ease of computation, this work assumes a uniform grid for fifth-order accuracy but also works on nonuniform grids, as will be shown in Section \ref{sec2d}. Due to the high order of accuracy, and the presented results living on (quasi-)uniform grids, this is not considered to have a large impact. Extension to arbitrary grids can be found in \citet{wellerAdaptivelyImplicitExplicitDensityWeightedsubmitted} and is flagged as a future generalisation here.

\subsection{Stability Analysis}\label{stability}

\begin{figure}[h!]
    \centering
    \includegraphics[width=0.49\textwidth]{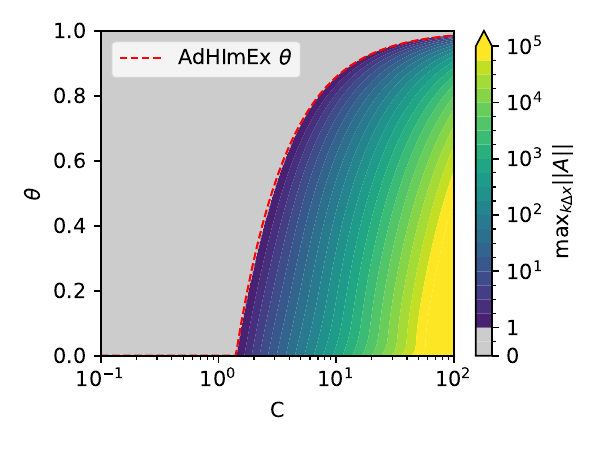}
    \caption{Stability of the AdHImEx scheme. Filled contours show the maximum absolute amplification factor for a discrete set of $C$ and $\theta$. Each maximum is taken over a range of $k\Delta x$. The AdHImEx $\theta$ (red dashes) lies in the grey stable region.}
    \label{figamplfac}
\end{figure}

In the above, we have been talking about the stability of AdHImEx. This is derived in this Section using von Neumann stability analysis. The idea is to find an amplification factor $A$ by assuming a 1D complex wave initial condition
\begin{equation}
    \Psi_i^{(n)} = A^{n} e^{\iota k x_i},
    \label{eqwave}
\end{equation}
where $\iota=\sqrt{-1}$, $k$ is the wavenumber, and $x_i=i\Delta x$. Applying this to a numerical scheme gives a criterion on $C$, as stability requires $||A(C)||\leq 1$. $A(k\Delta x, C, \theta)$ is found by substituting (\ref{eqbutcherAdHImEx}), (\ref{eqdiv}), (\ref{eqSD}), and (\ref{eqwave}) into (\ref{eqrungekutta}). Then for a pair of values of $C$ and $\theta$, assumed spatially uniform, we take a range of $k\Delta x$ values to find the maximum $||A(C,\theta)||$. Stability requires $||A(C,\theta)|| \leq 1$ for all $k \Delta x$: None of the wave modes are allowed to grow unstable. Hence, it suffices to focus on the maximum over $k\Delta x$. Figure \ref{figamplfac} shows the $\max_{k\Delta x}||A(C,\theta)||$ for the AdHImEx scheme for $C\in [0.1, 100]$ and $\theta \in [0,1]$. The maxima were taken over a discrete set of 401 $k\Delta x \in [0,2\pi]$. The combination of AdHImEx time stepping with fifth-order in space (\ref{eqSD}) is stable when
\begin{equation}
    \begin{aligned}
        \theta_c & = 1 - \frac{1}{1 + 0.7\max(0,C_c-1.4)},       \\
        \theta_f & = \max(\theta_{\text{u}}, \theta_{\text{d}}),
    \end{aligned}
    \label{eqtheta_adhimex}
\end{equation}
where
\begin{equation}
    C_c = \frac{1}{2} \frac{\Delta t}{V_c} \sum_{f\in c}|\vec{S}_f \cdot \vec{u}_f|.
    \label{eqcourant}
\end{equation}
The $c$ and $f$ subscripts refer to cell-centred and cell-face values, respectively. The $\text{d}$ and $\text{u}$ subscripts refer to the down- and upwind cells of the face, as in (\ref{eqSD}). The $1-1/C$ functional form emulates the implicitness for the first-order bounded AdImEx scheme \citep{wellerAdaptivelyImplicitMPDATA2023}. The $0.7$ and $1.4$ coefficients in (\ref{eqtheta_adhimex}) were found manually and checked numerically to satisfy $||A(C,\theta)||\leq 1$ for all discrete $|C| \leq 100$. (\ref{eqtheta_adhimex}) approaches $\theta=1$ asymptotically. Figure \ref{figamplfac} indicates $\theta$ as in (\ref{eqtheta_adhimex}) with a red dashed curve. It minimises $\theta$, thus maximising the use of the more accurate and efficient explicit time integration, while ensuring stability for large $C$. This stability was obtained based on cell-centred values, and the $\max$ in (\ref{eqtheta_adhimex}) ensures that stability is not broken when moving to face values when the Courant number is nonuniform in space. Stability under spatially varying $C$ and $\theta$ is only verified empirically (Section \ref{sec2d}).

\subsection{AdImEx Flux-Corrected Transport}\label{iterFCT}

A subset of the results in Sections \ref{sec1d} and \ref{sec2d} uses flux-corrected transport (FCT) to limit the solution fields, ensuring that the final result is either monotonic or nonnegative. This work uses a version of FCT that works for large time steps and nondivergent winds, described in \citet{wellerAdaptivelyImplicitMPDATA2023} and referred to as `AdImEx FCT'. AdImEx FCT gives a monotonic result, but has been found to be highly diffusive for large $C$ \citep[e.g.,][]{wellerAdaptivelyImplicitExplicitDensityWeightedsubmitted}. Hence, we also consider how the solution behaves when we solely ensure it is nonnegative. This is referred to as `NN FCT', and differs from AdImEx FCT only in terms of the locally allowable extrema: AdImEx FCT uses local values of the first-order upwind AdImEx solution at $(n+1)$ and possibly of $\Psi^{(n)}$ to determine these, whereas NN FCT uses global bounds $\Psi_{\min} = 0$ and $\Psi_{\max} \rightarrow \infty$.

\subsection{The WKS24 Scheme}

Results in the next Section compare AdHImEx to the AdImEx scheme from \citet{wellerAdaptivelyImplicitAdvection2024} (WKS24). Here, we review WKS24 applied to the advection equation. It has second-order accuracy in time for $C\leq2$, reducing to first-order accuracy for $C>2$, and has third-order accurate in space. It is a predictor-corrector scheme with two iterations $q=1,2$:
\begin{equation}
    \begin{aligned}
        \Psi^{q} = \Psi^{(n)} & - \Delta t \nabla \cdot \Big( (1-\alpha) \vec{u} \Psi^{(n)}_{\text{HO}}  \Big) - \Delta t \nabla \cdot \Big( \alpha (1-\beta) \vec{u} \Psi^{q-1}_{\text{1st}}  \Big) \\
                              & - \Delta t \nabla \cdot \Big( \alpha \beta \vec{u} \Psi^{q}_{\text{1st}} \Big) - \Delta t \nabla \cdot \Big( \alpha \vec{u} \Psi^{q-1}_{\text{HOC}}  \Big)
    \end{aligned}
\end{equation}
where $\Psi_{\text{HO}}  = \Psi_{\text{1st}} + \Psi_{\text{HOC}}$, $\alpha            = \max ( \frac{1}{2}, 1 - \frac{1}{C} )$, and $\beta             = 0 \text{ for } |C| < 0.8$ and 1 otherwise. Written in this form, $q=0$ and $2$ correspond to time levels $(n)$ and $(n+1)$ respectively.
This can be most readily compared to the original paper's (31)-(33), setting the Brunt-V\"ais\"al\"a frequency $N=0$ and the high-order flux limiter $\gamma=1$ as we are solving the linear advection equation. The quasi-cubic (HO) and first-order (1st) spatial discretisations are given similarly to (\ref{eqSD}) in 1D by
\begin{align}
    \Psi_{\text{HO}, f}  & = \tfrac{1}{3}\Psi_{\text{d}} + \tfrac{5}{6} \Psi_{\text{u}} - \tfrac{1}{6}\Psi_{\text{u}^2} \\
    \Psi_{\text{1st}, f} & = \Psi_{\text{u}},
\end{align}
along with (\ref{eqdiv}).

\subsection{The Matrix Solution}

The AdHImEx scheme requires the solution of a sparse linear system. The associated matrix can have nonzero values in seven diagonals: the main diagonal plus three neighbours on either side. The WKS24 scheme requires two sparse linear system solutions. For both schemes, the matrix elements are mostly zero due to the large explicit domain (where $\theta=0$).

Two methods are used to solve the matrix inversion. One-dimensional results solve the matrix inversion through the NumPy built-in solver, \texttt{numpy.linalg.solve}.
Two-dimensional results use the GMRES(m) iterative solver as described in \citet{saadGMRESGeneralizedMinimal1986}. This is run until convergence, minimising the residual to a solver tolerance relative to the Frobenius norm of the right-hand side of the matrix equation. This tolerance is $10^{-6}$ for AdHImEx and $10^{-15}$ for AdImEx FCT to provide monotonicity. The initial guess is $\Psi^{(n)}$ to ensure that constancy preservation is exact (Section \ref{res2dconstancy}). Each GMRES(m) restart takes four iterations. Furthermore, we require each matrix solution to take at least four iterations. This will be discussed in Section \ref{secruntime}.

Note that there is no matrix preconditioning and that the code has not been optimised for efficiency for this work. These are avenues of future optimisation. Matrix preconditioning is particularly important, as this could be shared between multiple tracers and thus support multi-tracer efficiency. The implicit solve itself does need to be performed for each tracer anew.

\section{One-dimensional Results}\label{sec1d}

In this Section, we show 1D AdHImEx results for uniform and nonuniform flows. Uniform flow is used for convergence analysis, comparison to WKS24, and the impact of both nonnegative and AdImEx FCT. The nonuniform flow provides a closer look at the removal of AdImEx artefacts through constancy preservation.

\subsection{Uniform Advection}\label{secuniform}

\begin{figure}[h!]
    \centering
    \includegraphics[width=\textwidth]{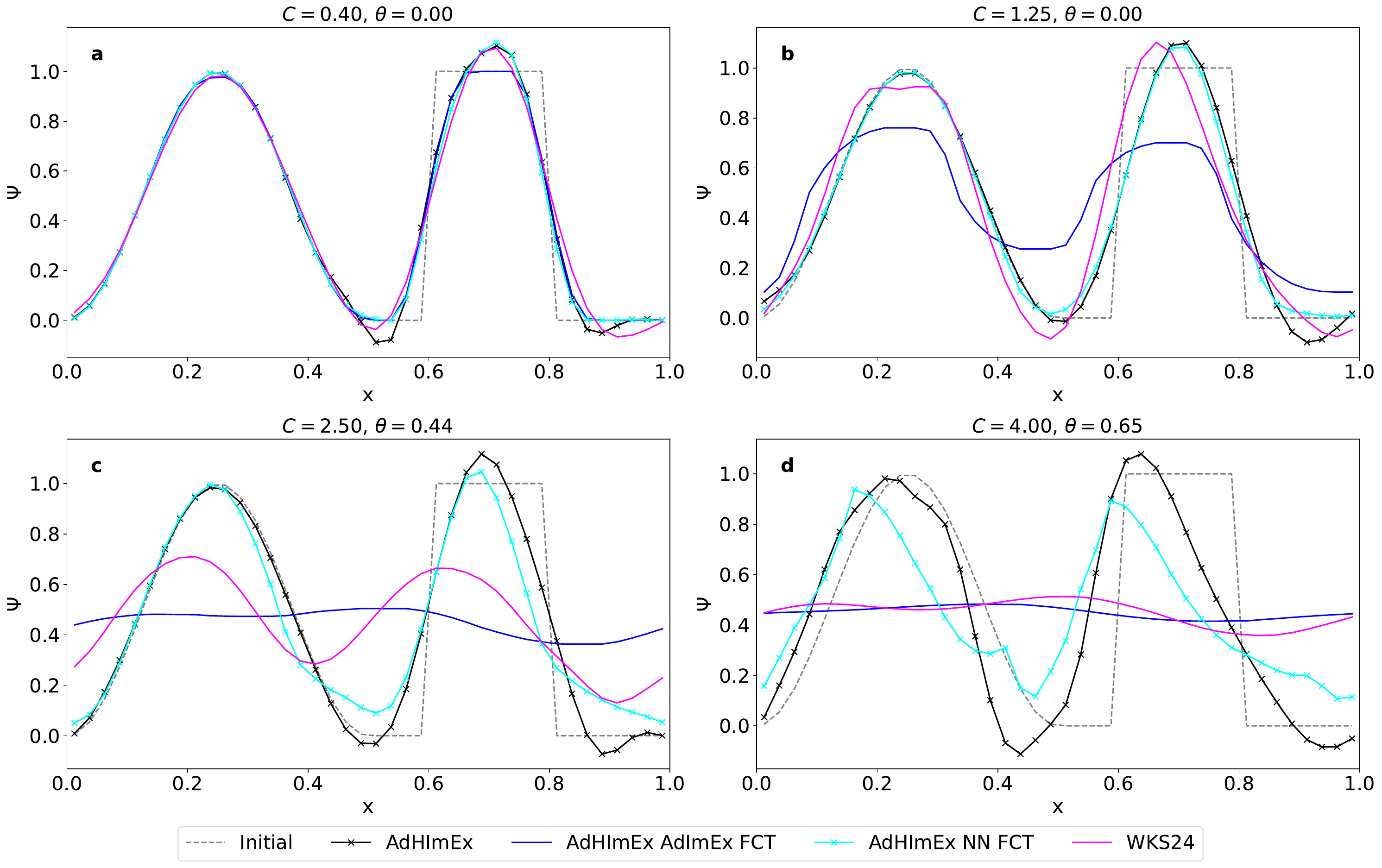}
    \caption{The skill of AdHImEx and the impact of FCT limiting as $C$ grows. $\Psi$ at $t=1$ with the cosine bell and square wave initial condition as in (\ref{eqICuniform}) in a periodic domain. Unlimited AdHImEx (black), AdImEx FCT- and NN FCT-limited AdHImEx, and the WKS24 scheme. $u=1$, $\Delta x=0.025$, and (a) $\Delta t=0.01$, $C=0.4$, (b) $\Delta t=0.03125$, $C=1.25$, (c) $\Delta t=0.0625$, $C=2.5$, and (d) $\Delta t=0.1$, $C=4$.}
    \label{uniform}
\end{figure}

The first set of results is intended to demonstrate the general quality and behaviour of the AdHImEx scheme and the impact of FCT limiting. Figure \ref{uniform} shows $\Psi$ after exactly one revolution over a periodic domain for various $C$. The initial condition is a combination of a cosine bell and a square wave of the domain, given by
\begin{equation}
    \Psi =
    \begin{cases}
        0.5(1-\cos (4\pi x)) & \text{for  } x \in [0,0.5],    \\
        1                    & \text{for  } x \in [0.6, 0.8], \\
        0                    & \text{otherwise.}
    \end{cases}
    \label{eqICuniform}
\end{equation}
The test is run with velocity $u=1$, grid spacing $\Delta x=$ 0.025, and number of grid points $n_x=$ 40. For the four subfigures, the number of time steps $n_t=$ 100, 32, 16, and 10 are varied along with $\Delta t =$ 0.01, 0.03125, 0.0625, and 0.1, resulting in $C=$ 0.4, 1.25, 2.5, and 4, respectively. The Courant numbers are constant in space and time. Figure \ref{uniform} shows results for the AdHImEx scheme without a limiter, and with the AdImEx and nonnegative FCT limiters.

As we see in Figure \ref{uniform}, AdHImEx without FCT is less diffusive and has smaller phase and amplitude errors for large $C$ than WKS24 without FCT. The AdImEx FCT-limited AdHImEx solution works well for small $C$, but gets increasingly diffusive for $C>1$. An important feature of FCT limiting is its independent treatment of in- and outfluxes. At each cell, FCT limits the in- and outfluxes separately, avoiding iterations to find the fluxes that provide the required limiting in all cells. However, this can lead to an overly limited result for AdImEx FCT for large $C$, since the balance between inflow and outflow is particularly important to avoid cells emptying or overflowing. As both the in- and outfluxes grow proportionally to $C$ while the allowable rise and fall do not change, FCT imposes increasingly strict limits on the high-order fluxes for growing $C$. This ensures no cell empties or overflows, but for AdImEx FCT with large $C$, the solution consequently largely follows the first-order bounded AdImEx scheme. In conclusion, AdHImEx itself works well for large $C$ but is let down by overly diffusive limiting. These results are consistent with existing literature: AdLImHEx, using first-order implicit time stepping, has a similar degradation for large $C$ \citep{wellerAdaptivelyImplicitExplicitDensityWeightedsubmitted}.

For NN FCT, the independent in/out flux limiting can amplify or create new extrema for large $C$. It might restrict a cell's large outward flux to avoid the creation of a negative value, but not touch its large inward flux, resulting in a new, spurious local maximum. However, apart from that, the NN FCT-limited AdHImEx follows the exact solution more closely and has smaller phase and amplitude errors than AdImEx FCT, even for $C>1$.

\subsection{Convergence} \label{seccourant_fixed}

\begin{figure}[h!]
    \centering
    \includegraphics[width=0.8\textwidth]{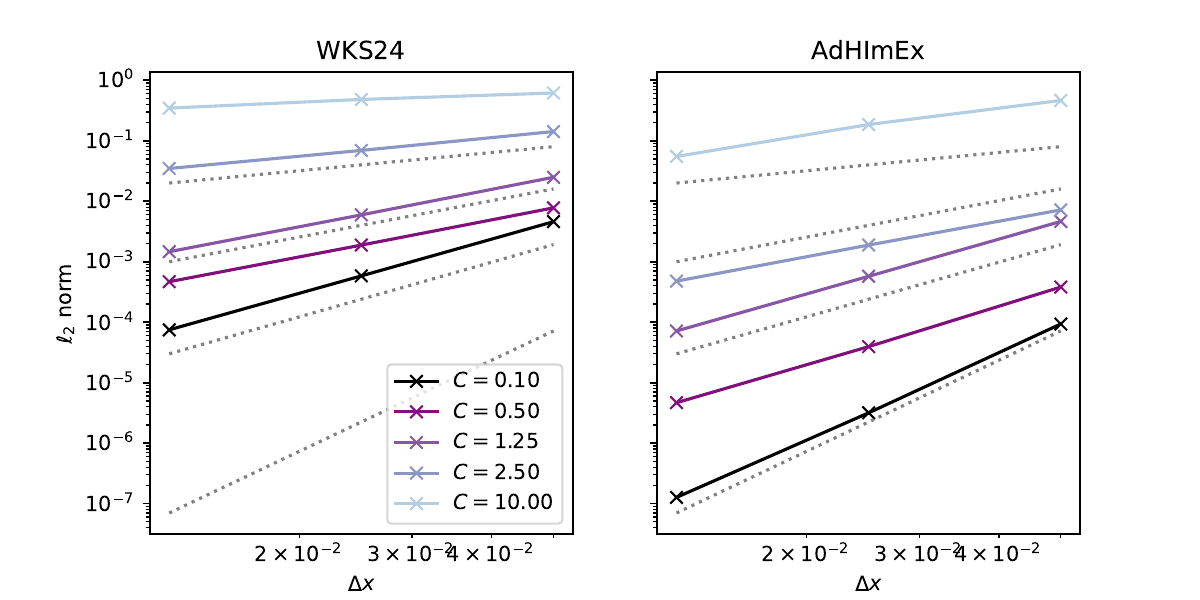}
    \caption{Order of convergence of AdHImEx (right) relative to WKS24 (left). $\ell_2$ error norms over $\Delta x$ for a sine wave initial condition, with $\Delta t$ and $n_t$ from Table \ref{tableconvergence}. First-, second-, third-, and fifth-order lines are dotted for reference.}
    \label{convergence}
\end{figure}

We present the order of convergence for various $C$ of AdHImEx and WKS24. Figure \ref{convergence} shows the $\ell_2$ error norms, defined here as root-mean-square errors, of numerical solutions in comparison to analytic solutions ($\Psi - \Psi_{\text{analytic}}$) over the whole domain at $t=0.5$. The analytic solution is the initial condition evaluated at $x-u(t-t_0)$, where the start time $t_0=0$ and the velocity $u=1$. The numerical solutions are calculated with $\Delta x=0.0125,0.025, 0.05$, and $n_x=80,40,20$ respectively, to give the three points on each line when combined with the $n_t$ and $\Delta t$ values in Table \ref{tableconvergence}. This gives five curves, with $C=$ 0.1, 0.5, 1.25, 2.5, 10. The initial condition is a sine wave: $\Psi = 0.5 + 0.5\sin (2\pi x)$ over an $x=[0,1)$ periodic domain.

For AdHImEx, as $C$ approaches zero, fifth-order spatial accuracy dominates the solution (Figure \ref{convergence}). As $C$ grows, the time discretisation starts to dominate the solution error: The fifth-order accuracy reduces to third-order, and further down to second-order for $C>1.4$ as expected. For WKS24, we see lower accuracy and order of convergence throughout: third-/second-order for small $C$, going down to first-order as $C>2$. For $C=2.5-10$, the AdHImEx solution is still significantly more accurate than WKS24, illustrating the utility of using second-order accurate time stepping.

\begin{table}[h!]
    \centering
    \begin{tabular}{| c | c ? c | c | c | c | c | c | c | c | c | c |}
        \hline
        \multicolumn{2}{|c}{} & \multicolumn{2}{c}{$C=0.1$} & \multicolumn{2}{c}{$C=0.5$} & \multicolumn{2}{c}{$C=1.25$} & \multicolumn{2}{c}{$C=2.5$} & \multicolumn{2}{c|}{$C=10$}                                                                \\
        \hline
        $\Delta x$            & $n_x$                       & $\Delta t$                  & $n_t$                        & $\Delta t$                  & $n_t$                       & $\Delta t$ & $n_t$ & $\Delta t$ & $n_t$ & $\Delta t$ & $n_t$ \\
        \thickhline
        0.0125                & 80                          & 0.00125                     & 400                          & 0.00625                     & 80                          & 0.015625   & 32    & 0.03125    & 16    & 0.125      & 4     \\
        0.025                 & 40                          & 0.0025                      & 200                          & 0.0125                      & 40                          & 0.03125    & 16    & 0.0625     & 8     & 0.25       & 2     \\
        0.05                  & 20                          & 0.005                       & 100                          & 0.025                       & 20                          & 0.0625     & 8     & 0.125      & 4     & 0.5        & 1     \\
        \hline
    \end{tabular}
    \caption{$\Delta t$ and $n_t$ values used in Figure \ref{convergence}, chosen to hold $C$ fixed along each convergence curve.}
    \label{tableconvergence}
\end{table}

\subsection{Impact of Constancy Preservation}\label{secartefacts}

An appropriate framework to address the removal of artefacts through constancy preservation in AdHImEx is given by 1D flow with varying Courant number. The numerical check of constancy preservation will follow in Section \ref{res2dconstancy}. Here, we consider its impact on artefacts by comparing the constancy-preserving version (using (\ref{eqFs})) of AdHImEx with its non-constancy-preserving counterpart (using (\ref{eqFs_ncp})). (\ref{eqFs}) was derived to preserve uniform fields over time by removing spurious divergence from Runge-Kutta stages 2 and 3 in the AdHImEx time stepping. We expect this to have an impact not only on uniform fields, but to remove artefacts in the solution more widely.

Consider a uniformly periodic domain with $x \in [0,1)$, $\Delta x = 0.004$, and $n_x = 250$. The initial condition and 1D nonuniform, divergent velocity are given by
\begin{equation}
    \begin{aligned}
        \Psi^{(n=0)} & = 50(1 + \sin(2 \pi (x - 0.3))), \\
        u            & = 0.2 + 9.8\sin^4(\pi(x + 0.5)),
    \end{aligned}
    \label{eqsubstages}
\end{equation}
resulting in a strongly yet smoothly spatially varying $C \in [0.2,10]$.

\begin{figure}[h!]
    \centering
    \includegraphics[width=\textwidth]{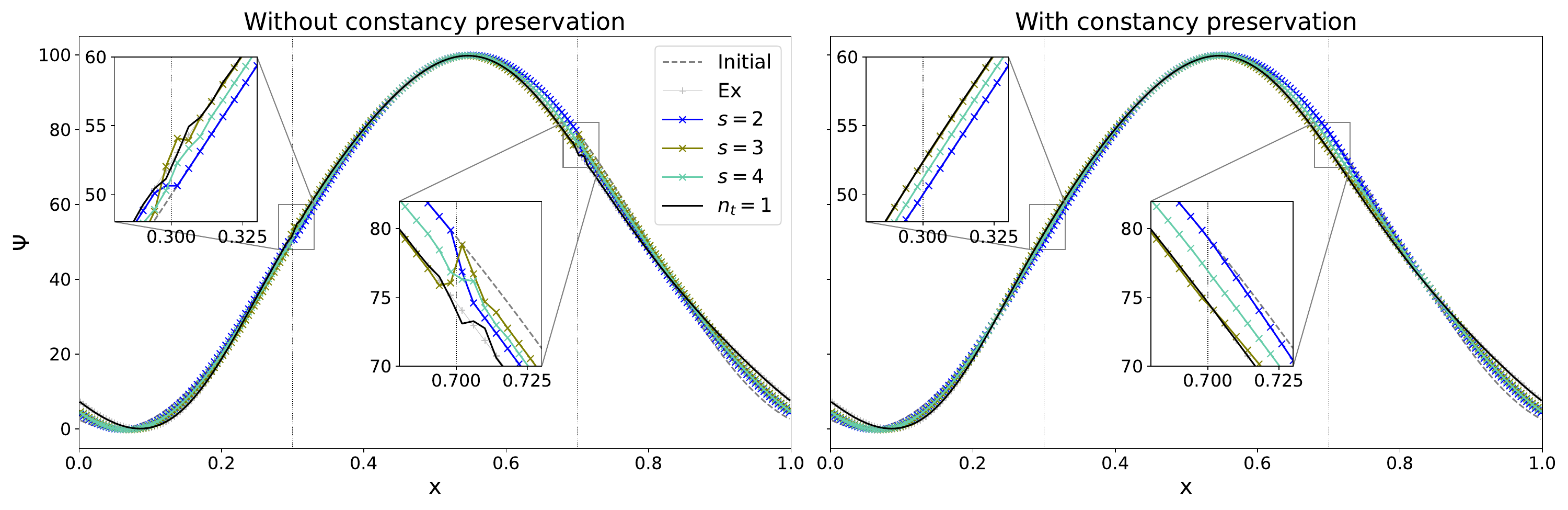}
    \caption{Removal of transport artefacts by constancy preservation. AdHImEx stage and final $\Psi$'s without (left, using (\ref{eqFs_ncp})) and with (right, using (\ref{eqFs})) constancy preservation, for one $\Delta t=0.004$ with a uniform $\Delta x=0.004$, and $\Psi^{(n=0)}$ and divergent $u$ from (\ref{eqsubstages}). The insets show the AdImEx interfaces at $x=0.3$ and $0.7$, where $\theta$ switches between zero and nonzero.}
    \label{substages}
\end{figure}

Figure \ref{substages} shows AdHImEx results after one time step $\Delta t=0.004$. The left and right Subfigures use (\ref{eqFs_ncp}) and (\ref{eqFs}), respectively. Both show the initial condition and explicit reference solution, alongside all AdHImEx stages. The explicit reference solution is third-order SSP Runge-Kutta with fifth-order in space, as defined by AdHImEx with $\theta=0$ globally. This is achieved by advancing the solution to $t=0.004$ using a time step reduced by a factor of 20. The $s$ stages coincide with the rows in the Butcher tableau in Section \ref{timestepping}. $s=1$ and 5 are not expressly labelled in Figure \ref{substages} as they are equivalent to the initial condition and $n_t=1$, respectively.

Using (\ref{eqFs}) removes artefacts for varying $\theta$, most notably around AdImEx interfaces $x=0.3$ and $0.7$ with a strong $\theta$ variation: $\theta$ switches between zero and nonzero. Without constancy preservation, artefacts from spurious divergence appear in stages 2 and 3 and are not completely removed by the end of the time step. With constancy preservation, these artefacts do not appear. This increases accuracy: The $\ell_2$ error norm (calculated with respect to the explicit reference solution at the same $t$) is about three times smaller with constancy preservation ($6.3 \times 10^{-4}$) than without ($1.9 \times 10^{-3}$).

\section{Two-dimensional Results}\label{sec2d}

The two-dimensional results consist of BK25 deformational flow \citep{bendallSWIFTMonotonicFluxForm2025} and Hadley-like circulation \citep{kentDynamicalCoreModel2014}. The BK25 nondivergent test allows us to consider the advection of smoothly and nonsmoothly varying tracers, addressing constancy preservation, accuracy, monotonicity, and efficiency. AdHImEx results are presented with both a uniform and nonuniform grid, the latter bringing about large Courant numbers through grid refinement. Hadley-like circulation lets us to look at the treatment of filaments in the solution and the skill of AdHImEx with a Courant number variation determined solely through variation in the wind. Note that all these results are on a Cartesian planar grid.

\subsection{BK25 Nondivergent Deformational Flow}\label{sec2dswiftnondiv}

We now consider the two-dimensional nondivergent deformational flow test case of \citet[][Appendix B]{bendallSWIFTMonotonicFluxForm2025}. This testcase allows us to assess the accuracy and boundedness of AdHImEx for nonsmooth tracer fields, while also providing a challenging environment for large Courant numbers on nonuniform meshes. We adopt a unity density, such that tracer concentration and tracer mass density are equivalent.

\begin{figure}[h!]
    \centering
    \includegraphics[width=\textwidth]{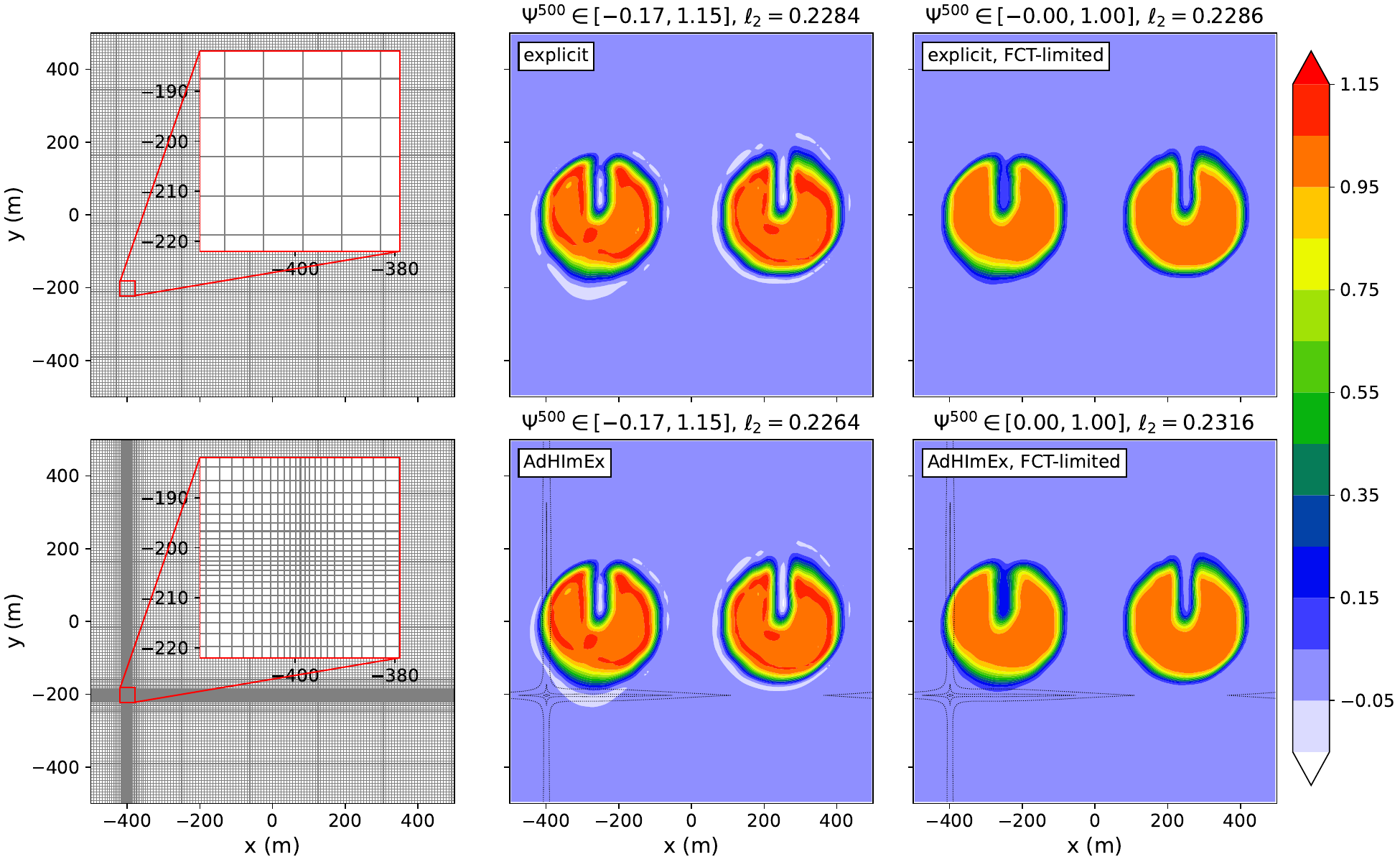}
    \caption{AdHImEx transport on a locally refined mesh with large Courant numbers, against an explicit reference. BK25 deformational flow for $\Psi$ at $t=100$ with $\Delta t=0.2$. The explicit reference (top row; AdHImEx with $\theta=0$) uses a uniform grid such that $C\in[0.18,0.84]\ \forall t$; the AdHImEx solution (bottom row) uses a nonuniform grid such that $C\in[0.18,6.34]\ \forall t$. These grids are shown in the left column, unlimited results in the middle column, and AdImEx FCT-limited results in the right column. Contours are shown at $C=[1.4,2.8,4.2,5.6]$.}
    \label{resswift2dnondiv}
\end{figure}

The tracer is advected in the positive $x$- and $y$-directions while undergoing deformation and returns to its initial condition after one period $T=100$. The test is run with $\Delta t=0.2$ and $n_t=500$ in a planar periodic domain with $(x,y)\in[-500,500]$\,m.

Results are presented on a $128\times128$ uniform grid and a $153\times153$ nonuniform grid. The uniform grid has $\Delta x=\Delta y=7.8125$\,m. The nonuniform grid has $\Delta x_{\max}=\Delta y_{\max}=7.8125$\,m and an idealised grid refinement down to $\Delta x_{\min}=\Delta y_{\min}=0.8696$\,m. This mimics the addition of grid points in regions of interest. For a purely explicit Eulerian advection scheme, such a local refinement would require a reduction in the time step for stability. With AdHImEx, this is not required.

The nonuniform grid matches the uniform grid throughout most of the domain, allowing differences between the two solutions to be attributed primarily to the refined region centred at $(x,y)\approx(-400,-200)$\,m. Grid compression is applied independently in each coordinate direction. In the $x$-direction (starting at $x=-500$\,m), there are 5 cells with a grid spacing of 7.8125\,m, followed by 20 cells that each multiply the grid spacing by 0.9, a smallest centre cell of 0.8696\,m, another 20 cells that each divide the grid spacing by 0.9, and finally 107 cells at 7.8125\,m. In the $y$-direction, the same compression is applied, but with 30 cells before the compression region and 82 afterwards. The grids are shown in the left column of Figure \ref{resswift2dnondiv}.

The nonuniform-grid experiment provides a stringent test of AdHImEx. Local mesh refinement simultaneously introduces large local Courant numbers and a spatial discretisation error, since the fifth-order stencil is formally derived for uniform meshes. Specifically, the spacing between neighbouring points from $\text{d}^2$ to $\text{u}^3$ in (\ref{eqSD}) is assumed constant. Consequently, deterioration in the solution may arise from both the time integration and the spatial discretisation. Despite these two challenges, the results presented below demonstrate that AdHImEx maintains a level of accuracy comparable to the fully explicit reference solution.

Figure \ref{resswift2dnondiv} shows results at $t=100$ for the slotted-cylinder initial condition described by BK25. The middle column shows unlimited AdHImEx and the right column shows AdHImEx with AdImEx FCT limiting. The top row shows explicit reference results on the uniform grid, for which $C\in[0.18,0.84]$. This reference is formed by AdHImEx with $\theta=0$ throughout. The bottom row shows results on the nonuniform grid, which requires AdHImEx time stepping and reaches $C\in[0.18,6.34]$. The largest Courant numbers occur in the smallest grid cells. Contours are shown at $C=[1.4,2.8,4.2,5.6]$.

Despite local Courant numbers exceeding six and the use of a nonuniform mesh, there is no visible deterioration in the AdHImEx solution compared with the fully explicit reference solution, either with or without AdImEx FCT limiting. AdImEx FCT removes newly created extrema, ensuring that the solution remains locally bounded. The $\ell_2$ error norm, comparing the tracer field at $t=100$ with its initial condition, remains essentially unchanged, even though the maximum Courant number increases from approximately $0.8$ to more than $6$. Without AdImEx FCT limiting, the nonuniform-grid result even exhibits a slightly smaller error than the uniform-grid reference because of its locally increased resolution. The excellent agreement with the explicit solution suggests that AdHImEx time stepping does not introduce an unexpected source of error. More broadly, the experiment shows that AdHImEx enables local mesh refinement without the time step restrictions normally associated with explicit Eulerian advection schemes.

\subsection{Constancy Preservation}\label{res2dconstancy}

Section \ref{secconstancy} discussed how AdHImEx keeps constant fields constant under nondivergent flow. BK25 nondivergent flow and the nonuniform grid from the previous Section can be employed to numerically confirm this. Instead of the slotted cylinders, we now initialise $\Psi$ with 0.5 everywhere. This $\Psi$ should not change as we step through time.

\begin{figure}[h!]
    \centering
    \includegraphics[width=0.7\textwidth]{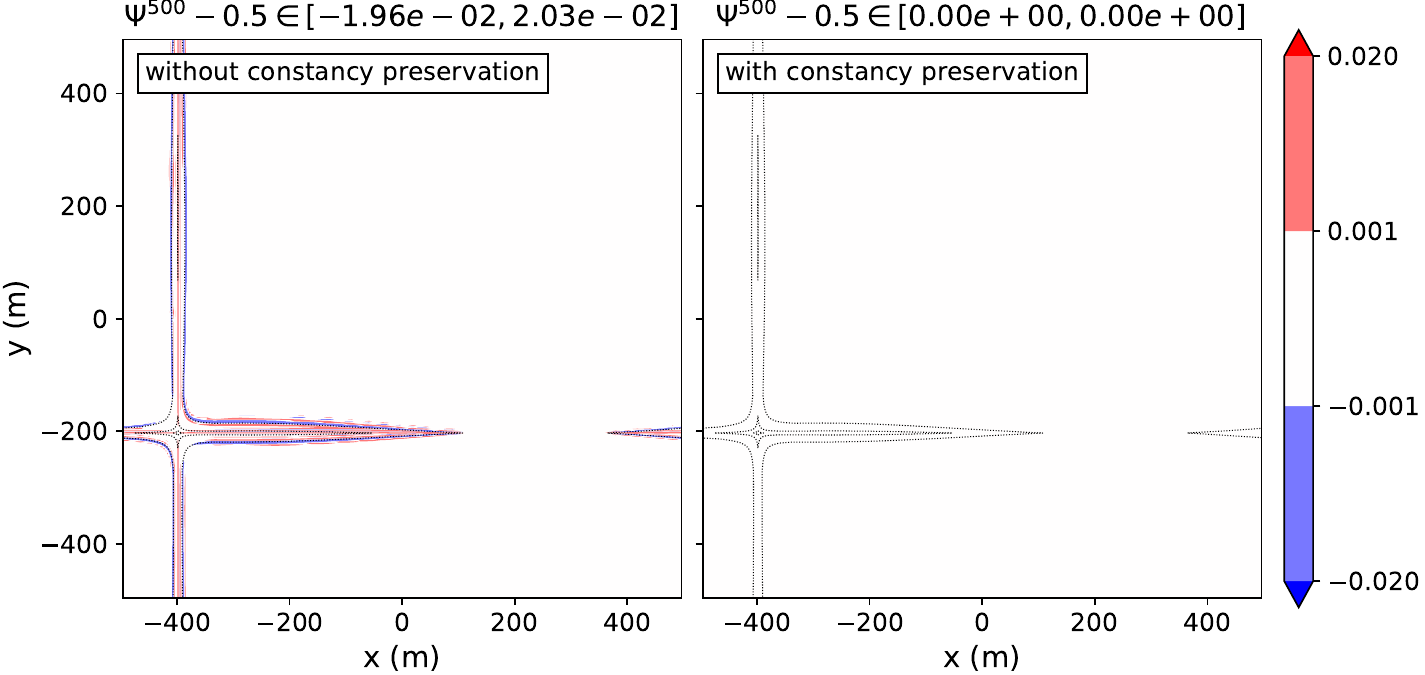}
    \caption{Exact numerical constancy preservation of a uniform field. $\Psi-0.5$ after advection by the BK25 deformational flow at $t=100$ with $\Delta t=0.2$ on the nonuniform grid from Section \ref{sec2dswiftnondiv}. Contours are shown at $C=[1.4,2.8,4.2,5.6]$. (left) AdHImEx without constancy preservation, using (\ref{eqFs_ncp}). Spurious divergence appears where $\theta$ varies in space. (right) AdHImEx with constancy preservation, using (\ref{eqFs}).}
    \label{resconstancy}
\end{figure}

Figure \ref{resconstancy} presents the result at $t=100$. Contours are shown at $C=[1.4, 2.8, 4.2, 5.6]$. The left Subfigure shows AdHImEx without constancy preservation, using (\ref{eqFs_ncp}). We see that the initially uniform field does not remain uniform whenever $\theta$ varies in space. The right Subfigure shows AdHImEx with constancy preservation, using (\ref{eqFs}). The result is exact because we use $\Psi^{(n)}$ as the matrix solution's initial guess.

\subsection{Efficiency}\label{secruntime}

\begin{figure}[h!]
    \centering
    \includegraphics[width=0.5\textwidth]{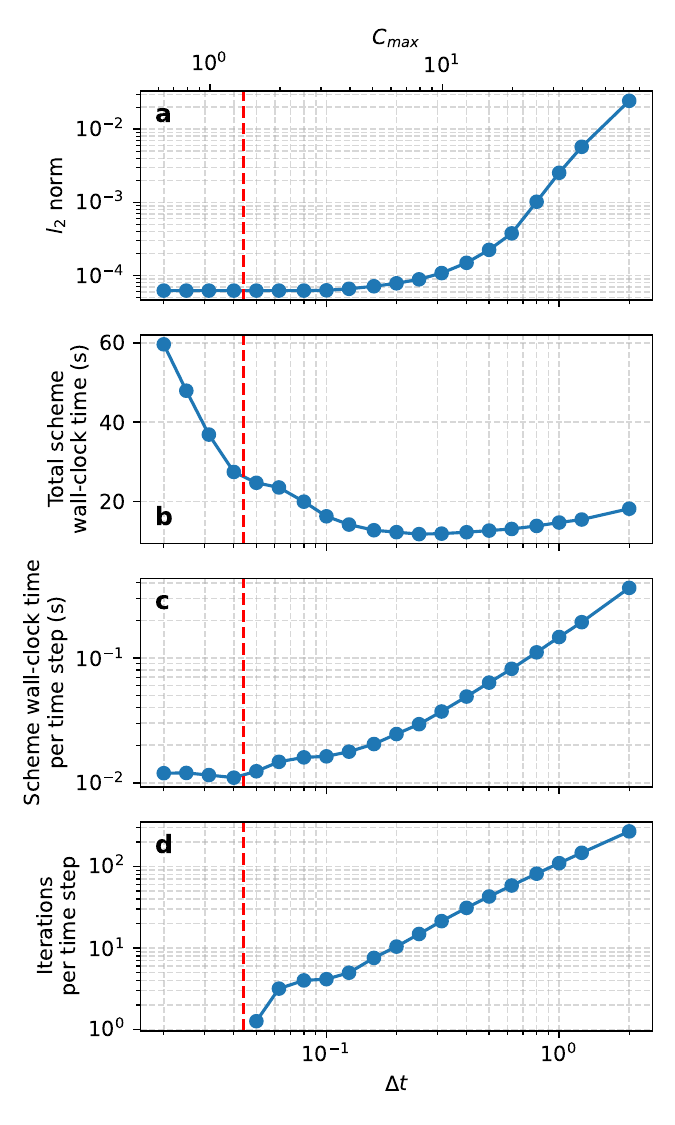}
    \caption{Computational cost and accuracy of AdHImEx as $\Delta t$ and $C$ increase. (a) $\ell_2$ error norms, (b) total wall-clock times, (c) wall-clock times per time step, and (d) iterations per time step for BK25 deformational flow and the nonuniform grid from Section \ref{sec2dswiftnondiv}. $n_t\in[5000,50]$ and $\Delta t\in[0.02,2]$\,s are varied to give $C_{\max}\in[0.634,63.4]$. Each point is an average of 20 simulations. The red dashes indicate $C_{\max}=1.4$, beyond which a matrix solution is required.}
    \label{restiming_accuracy}
\end{figure}

The efficiency of AdHImEx is assessed using the BK25 nondivergent flow setup from Section \ref{sec2dswiftnondiv} on the nonuniform grid. Unlike Section \ref{sec2dswiftnondiv}, which considered the transport of two slotted cylinders, the present experiment uses the double-sine-wave initial condition from BK25. Whereas Section \ref{sec2dswiftnondiv} considered only $n_t=500$, corresponding to $C_{\max}=6.34$, here $n_t\in [5000,50]$ and $\Delta t\in [0.02,2]$\,s are varied to obtain a range of $C_{\max}\in[0.634,63.4]$ while maintaining a fixed simulation length of $t=100$. Although computational cost is inherently case-dependent, this setup allows the efficiency of AdHImEx to be evaluated as a function of $C_{\max}$ and compared with fully explicit time stepping.

The central question is: Can it be more efficient to use AdHImEx rather than explicit time stepping? Results confirm that this is the case. Figure \ref{restiming_accuracy} shows the $\ell_2$ error norm (compared to the initial condition), AdHImEx's total wall-clock time, AdHImEx's wall-clock time per time step, and the total number of iterations taken throughout the simulation divided by the number of time steps. These are shown with respect to the time step size and maximum Courant number of the simulation. Each data point is an average of 20 identical simulations. The vertical dashed line indicates $C_{\max} > 1.4$, from which we start to have an implicit region, AdImEx time stepping, and thus a matrix solution. For the smallest $\Delta t$ (0.02\,s), the AdHImEx time integration is fully explicit and the advection scheme takes around 60\,s in total. As we increase $\Delta t$ to 2\,s, the total wall-clock time of the advection scheme is reduced to below 20\,s for $C_{\max} \gtrsim 3$, with the $\ell_2$ error norm remaining low (below $10^{-4}$) for $C_{\max} \lesssim 8$. The wall-clock time per time step increases, but this is more than balanced out by the smaller number of time steps taken. Notably, this does not take into account the reduced number of other updates required for, e.g., the velocity. These wall-clock times consider solely the advection scheme.

For the first few data points above $C_{\max}=1.4$, a relatively large number of iterations is observed because a minimum of four iterations is enforced for each matrix solution. Without this constraint, the accuracy of these cases deteriorates. The minimum of four iterations was found to provide a good compromise between computational cost and accuracy. Using a stricter convergence tolerance instead would unnecessarily increase the computational cost of simulations with larger $C_{\max}$.

These results use nonoptimised code and, in particular, the matrix is not preconditioned: Further improvements in favour of AdHImEx time stepping can reasonably be expected. For future work, it is recommended to combine the implicit time stepping with a third- or lower-order rather than a fifth-order spatial discretisation. This would increase the sparsity of the matrix and reduce the cost of the matrix solution. Even though these results already show efficiency gains when using AdHImEx time stepping, the fifth-order accuracy gains presumably do not balance out the increased computational cost of the solver. A better compromise might be achieved with a less-accurate, but cheaper implicit spatial discretisation. The associated sparser matrix would speed up the iterative solver and thus reduce the time-to-solution.

\citet{wellerAdaptivelyImplicitExplicitDensityWeightedsubmitted} set a high bar for efficiency of AdImEx time stepping schemes. The cost was nearly constant with $C$. For AdHImEx, this is not the case, but we do see that the cost per time step scales favourably (approximately logarithmically) with $C$. Besides that, various differences between the AdLImHEx and AdHImEx setups are relevant here: AdLImHEx uses a preconditioner, an optimised solver, and solves a much sparser and easier to solve matrix. As a result, it needs a smaller number of iterations to reach convergence. AdLImHEx does solve four matrices per time step, whereas AdHImEx solves only one.

\citet{wellerAdaptivelyImplicitAdvection2024} state that a common concern with implicit time stepping is the significant cost associated to the matrix solution. The expectation is that the cost increases linearly with the time step and thus there will be no efficiency gains. Here, we have shown that this is not the case for AdHImEx. In conjunction with the AdLImHEx efficiency results, it is clear that AdImEx time stepping can be competitive. Lastly, as mentioned in \citet{wellerAdaptivelyImplicitMPDATA2023}, although each tracer would require a matrix solve, they share the preconditioning step, aiding multi-tracer efficiency.

\subsection{Hadley-like Circulation}

Lastly, we consider a nondivergent Hadley-like circulation inspired by \citet{kentDynamicalCoreModel2014} to assess (1) the treatment of filaments and (2) the skill of AdHImEx time stepping without additional errors arising from the spatial discretisation. In the previous test cases, large Courant numbers were generated through variations in grid spacing. Consequently, whenever AdImEx time stepping was activated, the results reflected not only the effects of the temporal discretisation but also the loss of fifth-order spatial accuracy on a nonuniform grid. Here, the Courant number is varied solely through the velocity field, allowing the spatial discretisation to remain fifth-order accurate throughout the domain.

The testcase is adapted from \citet{kentDynamicalCoreModel2014} to a two-dimensional planar domain with latitude $\phi \in [-90^{\circ},90^{\circ}]$, height $z \in [0,z_{\text{top}}=12\,\mathrm{km}]$, and time $t \in [0,\tau=1\,\mathrm{day}]$.
The initial condition is
\begin{equation}
    \Psi^{(n=0)} = \begin{cases}
        0.5\Big( 1 + \cos \Big( \frac{2 \pi (z - z_0)}{z_2 - z_1} \Big) \Big) & \text{  if  } z_1 \leq z \leq z_2 \\
        0                                                                     & \text{  otherwise.  }
    \end{cases}
\end{equation}
where $z_0 = 0.5(z_1 + z_2)$, $z_1 = 2 \text{ km}$, and $z_2 = 5 \text{ km}$. The streamfunction is given by
\begin{equation}
    \psi = \frac{w_0}{K} \cos (\phi) \sin ( K \phi ) \sin \Big( \frac{\pi z}{z_{\text{top}}} \Big) \cos \Big( \frac{\pi t}{\tau} \Big)
\end{equation}
with $K=5$ and $w_0 = 0.15$\,m/s.
The original setup of \citet{kentDynamicalCoreModel2014} is three-dimensional with a longitudinally uniform flow, making it straightforward to adapt to a two-dimensional pole-to-pole slice. The meridional and vertical winds are
\begin{equation}
    v = -\frac{\partial \psi}{\partial z}, \quad w = \frac{\partial \psi}{\partial \phi}.
\end{equation}
There is no flux through the boundaries. Ghost cells, with nearest-grid-point physical values, are used for the spatial discretisation at the boundaries. Unlike \citet{kentDynamicalCoreModel2014}, a constant unity air density is assumed throughout the domain. Furthermore, the meridional volume factor associated with spherical geometry is neglected, resulting in a latitude-independent meridional grid spacing. Consequently, this testcase should be regarded as a simplified Hadley-like circulation inspired by \citet{kentDynamicalCoreModel2014} rather than an exact reproduction. The constant-density assumption in particular modifies the deformation of the tracer field, leading to some qualitative differences from the solutions presented by \citet{kentDynamicalCoreModel2014}.

\begin{figure}[h!]
    \centering
    \includegraphics[width=\textwidth]{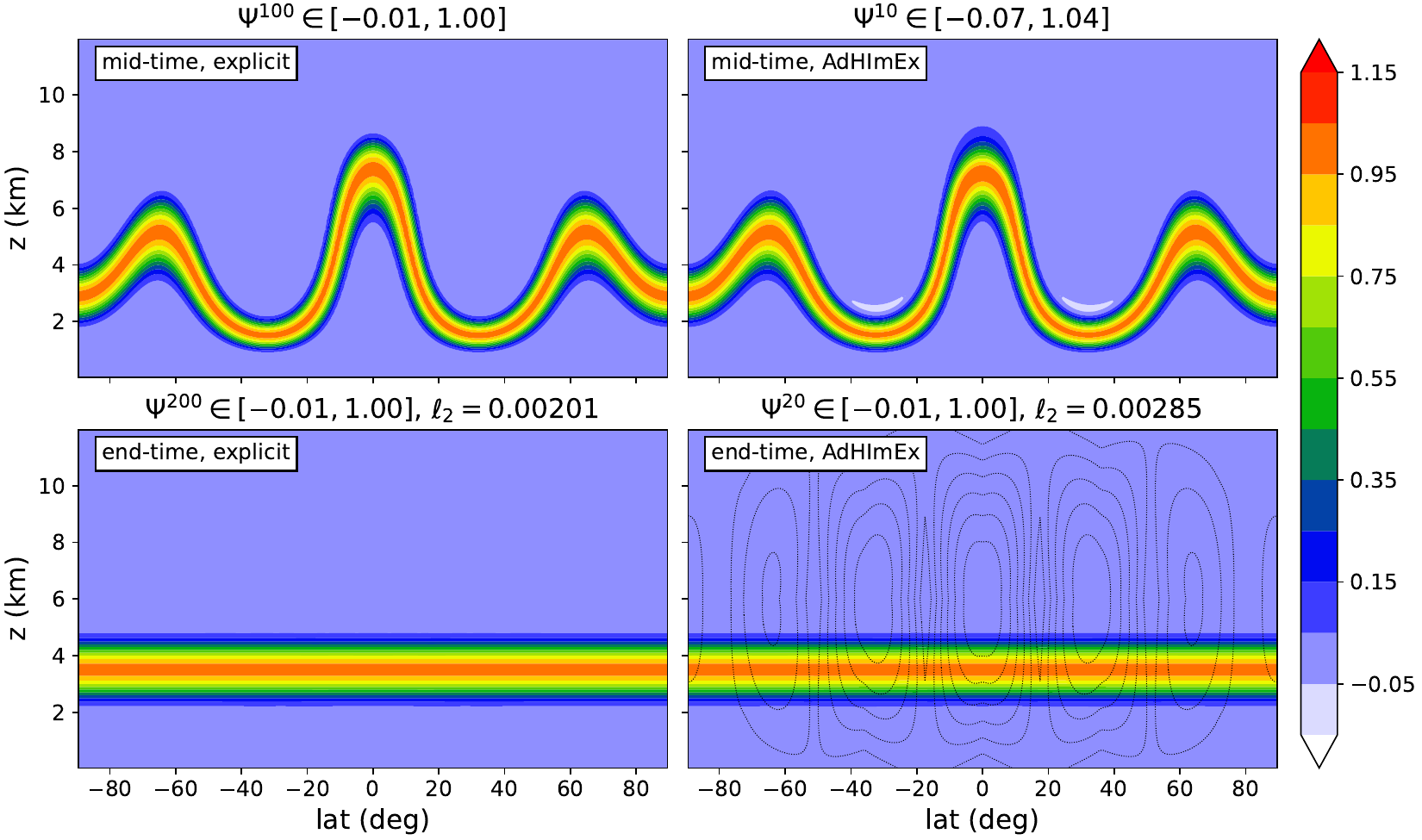}
    \caption{The skill of AdHImEx under strong filament deformation, with the Courant number varied solely through the winds. Hadley-like circulation for (left) explicit reference, $\Delta t=432$\,s such that $C\in[0.00,0.97]$, (right) AdHImEx with $\Delta t=4320$\,s such that $C\in[0.00,9.68]$. (top) Mid-time result at maximum deformation, and (bottom) end-time result ($t=1$\,day) with a return to the initial condition. Contours are shown at $C=[1.4,2.8,4.2,5.6,7,8.4]$.}
    \label{reshadley}
\end{figure}

This test is run with $n_\phi = n_z = 180$ for one day. Figure \ref{reshadley} compares the AdHImEx solution (right) to its explicit reference (left). The AdHImEx solution uses $n_t=20$ and $\Delta t = 4320\,$s such that $C\in[0.00,9.68]$. The explicit reference is AdHImEx with $\theta=0$ everywhere. It uses $n_t=200$ and a reduced time step of $\Delta t = 432\,$s such that $C\in[0.00, 0.97]$. The result is shown halfway through the simulation, at maximum deformation ($n_t=100$ for explicit, $n_t=10$ for AdHImEx), and at the end ($n_t=200$ and $n_t=20$), where there is a return to the initial condition. Contours are shown at $C=[1.4, 2.8, 4.2, 5.6, 7, 8.4]$. The mid-time AdHImEx solution exceeds the original bounds more than the explicit reference, but this largely reverts at the end of the testcase. This could be mitigated with limiting. The final AdHImEx and explicit results match each other closely. The increase in $\ell_2$ error norm (compared to the initial condition) for AdHImEx compared to the explicit reference is attributed to (1) the larger time step used for the AdHImEx solution and (2) the large part of the domain that has nonzero $\theta$ at the start and end of the simulation, as illustrated by the $C$ contour lines at the end-time. In conclusion, when fifth-order spatial accuracy is maintained throughout the domain, AdHImEx reproduces the explicit solution closely despite Courant numbers approaching 10. This indicates that the method introduces no substantial additional errors beyond those expected from the larger time step, even under strong filament deformation.

\section{Summary and Conclusions}\label{sec:conclusion}

This work set out to address two limitations of existing AdImEx schemes: first-order temporal accuracy in the implicit regime and the need for multiple implicit solves per time step. AdHImEx addresses both. The scheme has second-order temporal accuracy for large Courant numbers while requiring only a single implicit solve per time step, and is combined here with a conservative fifth-order finite-volume discretisation.

A second contribution is a stage-dependent treatment of the divergence operator that preserves both conservation and constancy when the implicitness varies in space and the implicit and explicit time stepping use different stage time steps. It avoids spurious divergence and associated transport artefacts. Crucially, this approach is not restricted to AdHImEx but is applicable to other AdImEx schemes.

The numerical results demonstrate substantial improvements over previous first-order AdImEx schemes. One-dimensional results confirmed the expected third- and second-order behaviour and significantly improved accuracy over WKS24 at large Courant numbers, with close agreement with explicit reference solutions.

In two dimensions, AdHImEx maintained explicit-reference accuracy despite large local Courant numbers. The BK25 nonuniform-grid experiment, in particular, demonstrates that local mesh refinement can be introduced without the restrictive time step reductions normally associated with explicit Eulerian transport.

Efficiency experiments showed that, although the cost per time step increases as the implicit region grows, the reduction in the number of time steps can lead to an overall decrease in computational cost. These results indicate that AdHImEx can provide efficiency gains compared to explicit transport schemes.

AdImEx FCT was found to become increasingly inaccurate for large Courant numbers, causing the solution to approach the underlying first-order bounded scheme. NN FCT retains much more of the high-order accuracy but may introduce artefacts. Limiting for large Courant numbers therefore remains an important area for future development.

Relative to AdLImHEx, AdHImEx improves the large-$C$ temporal accuracy from first- to second-order and reduces the number of implicit solves from four per time step to one. These gains come at the cost of a more expensive matrix solution. Combining the favourable accuracy properties of AdHImEx with the solver performance of AdLImHEx would be a valuable direction for future research.

Another clear gap of the present work is its restriction to one- and two-dimensional planar tests. Extension to spherical geometries, arbitrary meshes, density-weighted transport, and further optimisation of the solver remain as future work.

\section*{Acknowledgements}

This work was supported by the University of Reading Advancing the Frontiers of Earth System Prediction (AFESP) Programme [award number A3720300].

\section*{Code Availability}

Figures (\ref{figamplfac})-(\ref{substages}) are produced by running \texttt{python main.py} from \url{https://github.com/ambertewinkel/AdHImEx/tree/paper\_figures}. Figures (\ref{resswift2dnondiv})-(\ref{reshadley}) and the graphical abstract are produced by running \texttt{run\_paper\_results.sh} from \url{https://github.com/ambertewinkel/2Dadvection/tree/paper\_figures}.

\section*{Declaration of generative AI and AI-assisted technologies in the manuscript preparation process}

During the preparation of this work, the authors used Microsoft 365 Copilot and Claude for sentence-level text editing and figure formatting. After this, the authors reviewed and edited the content as needed and take full responsibility for the content of the published article.

\bibliography{mylibrary}

\end{document}